\documentclass{amsart}

\usepackage{amsmath}
\usepackage{hyperref}
\usepackage{amssymb}

\newtheorem*{theorem*}{Theorem}

\newcommand{\C}{\mathbb{C}}

\newcommand{\R}{\mathbb{R}}

\DeclareMathOperator{\Hom}{Hom}
\DeclareMathOperator{\Gal}{Gal}
\DeclareMathOperator{\im}{im}

\begin{document}

\title{Grothendieck's use of equality}
\author{Kevin Buzzard}
\email{k.buzzard@imperial.ac.uk}
\address{Department of Mathematics, Imperial College London}
\begin{abstract}
  We discuss how the concept of equality is used by mathematicians (including Grothendieck), and what effect this has when trying to formalise mathematics. We challenge various reasonable-sounding slogans about equality.

\end{abstract}

\maketitle

\section*{Overview}

Many mathematical objects and constructions are \emph{uniquely characterised} by some kind of defining property. For example the real numbers are the unique complete Archimedean ordered field, and many constructions in algebra (localisations, tensor products,\ldots) and topology (product of topological spaces, completion of a metric or uniform space,\ldots) are characterised uniquely by a universal property. To be more precise, I would like to discuss properties which define a mathematical object \emph{up to unique isomorphism}. This is a very strong statement: to give a non-example, we could consider the property of being a group of order 5. There is only one such group, up to isomorphism, however this group has automorphisms (for example the map sending an element to its square), meaning that it is unique but only up to non-unique isomorphism, and hence the property of being a group of order 5 is not the kind of property that this paper is about.

Let us now fix a property $P$ which uniquely characterises a mathematical object, up to unique isomorphism. Here are three assertions.

\smallskip

{\bf Assertion 1.} You can (and probably should) develop the theory of objects characterised by $P$ using only the assumption $P$, and you do not need to worry about any particular construction of the object (beyond knowing that \emph{some} construction exists which satisfies $P$, and hence that your theory is not vacuous).

\smallskip

{\bf Assertion 2.} Any two mathematical objects satisfying $P$ may in practice be assumed to be \emph{equal}, because any mathematically meaningful assertion satisfied by one is also satisfied by the other.

\smallskip

{\bf Assertion 3.} More generally, if you have two objects which are \emph{canonically} isomorphic, then these objects may in practice be assumed to be equal.

\smallskip

Let us consider an example. Whilst there are several different constructions of the real numbers from the rationals (Cauchy sequences, Dedekind cuts, Bourbaki uniform space completion,\ldots), all of the basic analysis taught in a first undergraduate course is developed using only the Archimedean and completeness axioms satisfied by the reals. A mathematician would not dream of saying \emph{which} definition of the real numbers they were using -- this would be absurd. Indeed Newton, Euler and Gauss were happily using ``the'' real numbers long before Cauchy, Dedekind and Bourbaki came along with their different models: each of their definitions of the ordered field $\R$ is \emph{uniquely isomophic} to the others, so which one we are \emph{actually} using doesn't matter in practice.

I will say more about localisations later, but there are also several different constructions of the localisation $R[1/S]$ of a commutative ring at a multiplicative subset (a quotient of $R\times S$, a quotient of a multivariable polynomial ring\ldots) and a mathematician would never state precisely which construction they are using when they write $R[1/S]$; we know that this is irrelevant.

In this paper I argue that the first assertion above is false, the second is dangerous, and the third is meaningless. Assertions 2 and 3 seem to be used in many places in the algebraic geometry literature -- indeed we will discuss Grothendieck's usage of the $=$ symbol in some depth later.

However, as a working mathematician I am also aware that there is something deeper going on here, which I find it difficult to put my finger on. Even though the word ``canonical'' has no formal meeting, mathematicians are certainly not being mindless in their use of the idea. The \emph{concept} that two objects are canonically isomorphic and can hence be \emph{identified} is an extremely important one in practice; it is a useful organisational principle, it reduces cognitive load, and it does not in practice introduce errors in arguments, at least if the author knows what they are doing. Similarly, making definitions and proving theorems about objects by using an explicit construction rather than the universal property is a practical tool which is used across mathematics (for example, picking a basis of a vector space to make a definition and then observing that the definition is independent of the choice).

What is frustrating about the situation is the following. Let's say that one is attempting to formalise some mathematics on a computer (that is, translating the mathematics from the paper literature into the language of an interactive theorem prover -- a computer program which knows the axioms of mathematics and the rules of logic). Right now this process involves writing down many of the details of one's arguments. When one is forced to write down what one actually \emph{means} and cannot hide behind such ill-defined words as ``canonical'', or claims that unequal things are equal, one sometimes finds that one has to do extra work, or even rethink how certain ideas should be presented. Indeed, sometimes the most painless way to do this work involves having to create new mathematical ideas which are not present in the paper arguments. I am certainly \emph{not} arguing that the literature is incorrect, but I am arguing that many arguments in the literature are often not strictly speaking complete from a formalist point of view. With the advent of the formalisation of the mathematics around algebraic and arithmetic geometry using computer theorem provers, for example the work coming out of the Lean prover community (\cite{bcmperfectoid}, \cite{livingstongpcoh}, \cite{defrutosfernandeznorms}, \cite{angdinata_yu}, \cite{zhang_proj},\ldots) and related work in Isabelle (\cite{isabelle-schemes}) and cubical Agda (\cite{zeuner2023univalent}), these things will start to matter. In section~\ref{five} below I give an explicit example of some real trouble which we ran into when proving a very basic theorem about schemes, and which was only ultimately resolved in a satisfactory way after some new mathematical ideas were developed by the Lean community.

\section{Acknowledgements}

This paper has its origins in a talk~\cite{chapmanvid} which I gave at Chapman University at the ``Grothendieck, a Multifarious Giant'' conference in 2022. Much of the paper was written in Coffee Zee, a coffee shop on the Holloway Road; I thank Seb for his excellent decaf oat cappucinos (even though I know he doesn't approve). Many thanks to both Brian Conrad and the referee for their comments on a preliminary version. All remaining errors are solely the fault of the author.

\section{Introduction}

Six years ago, I thought I understood mathematical equality. I thought that it was one well-defined term, and that there was nothing which could be said about it which was of any interest to me as a working mathematician with a knowledge of, but no real interest in, the foundations of my subject. Then I started to try and do masters level mathematics in a computer theorem prover, and I discovered that equality was a rather thornier concept than I had appreciated. In particular I discovered that computer scientists had, many years ago, isolated several different concepts of equality, and had a profound understanding of the subject. The three-character string ``$2+2$'', typed into a computer algebra system, is not equal to the one-character string ``$4$'' output by the system, for example; some sort of ``processing'' has taken place. A computer scientist might say that whilst the \emph{numbers} $2+2$ and $4$ are equal, the \emph{terms} are not. Mathematicians on the other hand are extremely good at internalising the processing and, after a while, ignoring it. In practice we use the concept of equality rather loosely, relying on some kind of profound intuition rather than the logical framework which some of us believe that we are actually working within. In fact we are far ahead of the computer scientists in these matters: the concept of mathematical equality (and in particular its usage in algebraic geometry to represent canonical isomorphism) is an idea which has not yet been captured by any of the formal definitions of the $=$ symbol that we see in the literature or in computer theorem provers. The set-theoretic definition is too strong (Cauchy reals and Dedekind reals are certainly not equal as sets) and the homotopy type theoretic definition is too weak (things can be equal in more than one way, in contrast to Grothendieck's usage of the term). Mathematicians might claim to have captured some kind of ``sweet spot'' -- but in fact, what I now believe is happening, is that we are using notational tricks to sweep various issues under the carpet rather than developing a more conceptual framework which would deal with them in a manner much more amenable to formalisation.

Fast forward to the present and I am still trying to do mathematics in a computer theorem prover (the Lean interactive theorem prover), rather than on pen and paper. I do this for the same reason that I write my papers in LaTeX rather than using a typewriter: I believe that it represents progress. Unfortunately, at the time of writing, the libraries of mathematics built in the various interactive theorem provers which currently exist are woefully inadequate for doing most modern mathematics. However, over the last few years myself and hundreds of other mathematicians in the Lean community have spent many many thousands of person-hours building a digitised library {\tt mathlib} (\cite{mathlib}) of standard undergraduate, Masters and early PhD level mathematics, so this is going to change. I hope that before I die, these computer tools will have matured to the extent that it is as easy to do mathematics in them as it is to currently do it on paper. Once the systems are this mature, it might be the case that future generations of mathematicians will not have to worry about what people like Grothendieck mean by equality, because the systems will allow us to use the concept the way that it is currently used by mathematicians in practice. They will also point out to the author cases when it turns out that they didn't fully perceive what they were doing (by pointing out possibly nontrivial issues which need to be resolved in order to make an argument complete).

\section{Universal properties.}

I am aware of three reasonable classes of foundations for mathematics: those based on set theory, those based on type theory, and those based on category theory. There are also various results saying that, broadly speaking, these foundations are able to prove the same theorems. However for a discussion of equality to take place we will have to pin down exactly what we are talking about, and so I will choose ZFC set theory and classical logic as the basis for our discussion.\footnote{Personally I now think that dependent type theory is a far more convenient foundation than set theory for mathematics, at least when it comes to doing mathematics in a computer theorem prover. For example, it greatly decreases the amount of ``junk questions'' which one can ask -- there is no type-theoretic analogue of the set-theoretic question ``is $0\in 1$?'', which (using the standard set-theoretic definitions) is true if $0$ and $1$ are natural numbers, but false if they are real numbers. In type theory this question cannot even be \emph{stated}, as the concept of ``an element'' and ``a set'' are different things (a term, and a type). Some believe that AI tools are one day going to become useful for human mathematicians (and eventually surpass them, although currently this is science fiction). Restricting junk theorems is just one way of pointing these tools in the right direction.} The basis for this decision is simply that if a mathematician ever went to a class on the logical foundations of their subject (and many of them did not), then it was likely to have been a set theory class. Moreover, the pure mathematics literature is written in a superficially set-theoretic style: we are told that group is a set with some structure and axioms, a manifold is a set with a different structure and axioms, and so on. Here the word ``set'' is just a placeholder for the idea of a ``collection of atoms''.

It is obvious that two sets which are equal have the same elements; this follows from the so-called \emph{principle of substitution} for equality, which states that if $X$ and $Y$ are any two mathematical objects which are equal, then any claim which you can make in your foundational system about $X$ is also true for $Y$. The converse, that two sets with the same elements are equal, is imposed as an axiom of the theory. This ensures that the abstract concept of a set coincides with our mental model of what it is representing: a set is no more and no less than a collection of stuff.

Now I would like to begin the discussion of properties which uniquely characterise a mathematical object. Let us start with an example: the \emph{product} $X\times Y$ of two sets~$X$ and~$Y$. Let me warn the reader now that in the next few paragraphs I will be making a very careful distinction between the concept of \emph{the product} of~$X$ and~$Y$, and the concept of \emph{a product} of~$X$ and~$Y$.

\emph{The} product $X\times Y$ of two sets is defined to be the set of ordered pairs $(x,y)$ with $x\in X$ and $y\in Y$ (one can check using the axioms of set theory that it is possible to create this set). Note that here we run into the same issue that we saw earlier with the real numbers: there are several distinct ways to define the concept of an ordered pair in set theory. Set theory is designed very well to work with unordered pairs: the sets $\{x,y\}$ and $\{y,x\}$ have the same elements and are thus equal, so to define an ordered pair one needs to use some kind of hack. The Wikipedia page for ordered pairs~\cite{wikipedia_ordered_pair} currently gives three distinct constructions, due to Wiener ($\{\{\{x\},\emptyset\},\{\{y\}\}\}$), Hausdorff ($\{\{x,1\},\{y,2\}\}$) and Kurotowski ($\{\{x\},\{x,y\}\}$); all have the air of being slightly contrived.\footnote{One might also attempt to define an ordered pair as a function $f:\{1,2\}\to X\cup Y$ with $f(1)\in X$ and $f(2)\in Y$. However this definition is circular if we use the usual conventions of ZFC set theory, because a function is defined to be a set of ordered pairs.} Again mathematicians are well aware that this issue does not matter at all in practice: all that we need to know is the defining property of ordered pairs, which is that $(x_1,y_1)=(x_2,y_2)\iff x_1=x_2$ and $y_1=y_2$; this is all that we shall need, and all the definitions satisfy this property.

The product $X\times Y$ is equipped with two projection maps $\pi_1:X\times Y\to X$ and $\pi_2: X\times Y\to Y$. Stricly speaking, it is not just the product, but the triple $(X\times Y, \pi_1, \pi_2)$ which satisfies the following \emph{universal property}:

\medskip

{\bf The universal property of products:} A triple $(P,\pi_1:P\to X,\pi_2:P\to Y)$ is called \emph{a product of $X$ and $Y$} if it satisfies the following property: if $S$ is any set at all, and $f:S\to X$ and $g:S\to Y$ are functions, then there is a unique function from $S$ to $P$ such that its composition with~$\pi_1$ is~$f$ and its composition with~$\pi_2$ is~$g$.

\medskip

The universal property is not a definition of \emph{the product} of two sets; it can be thought of as infinitely many facts which \emph{a product} needs to satisfy (one for each choice of set $S$ and functions $f$ and $g$). It is not hard to verify that \emph{the product} $X\times Y$ of $X$ and $Y$, equipped with the natural projections, is \emph{a product}. But the converse is not at all true: there are typically plenty of other triples $(P,\pi_1,\pi_2)$ which satisfy the property of being \emph{a product} without being \emph{the product}. For example, if $X=\{37\}$ and $Y=\{42\}$ then \emph{the product} $X\times Y$ is $\{(37,42)\}$, but in fact \emph{any} set $P$ with one element, equipped with $\pi_1:P\to X$ sending everything to 37 and $\pi_2:P\to Y$ sending everything to 42, satisfies the universal property of being \emph{a product}. In particular, there are in general uncountably many \emph{different} things which satisfy the property of being \emph{a product}.

However mathematicians are extremely good at \emph{identifying} these different things; they are ``the same'' in a manner which \emph{transcends the correct usage of the $=$ symbol}. We can do because of the \emph{uniqueness yoga for universal objects}. Let us go through this yoga, which is a piece of formal category-theoretic nonsense, in the case of products.

Say $P_1$ and $P_2$ are both a product for $X$ and $Y$. Applying the existence part of the universal property for $P_2$ (with its projections to $X$ and $Y$) to the set $S=P_1$ (with its projections to $X$ and $Y$) gives us a function $\alpha:P_1\to P_2$ commuting with the projections to $X$ and $Y$. Switching the $1$s and $2$s in the argument we can also construct a function $\beta:P_2\to P_1$ commuting with the projections. Moreover, $\beta\circ\alpha$ is a map from $P_1$ to $P_1$ commuting with the projections, as is the identity function; by the uniqueness part of the universal property of $P_1$ applied to $P_1$ we deduce that $\beta\circ\alpha$ must be the identity map on $P_1$; similarly, $\alpha\circ\beta$ must be the identity map on $P_2$. Hence $\alpha$ and $\beta$ are bijections commuting with the projection maps. Finally, using the uniqueness part of the universal property for $P_2$ applied to $P_1$ tells us that $\alpha$ is the unique map from $P_1$ to $P_2$ which commutes with the projections, by symmetry $\beta$ is the unique map from $P_2$ to $P_1$ which commutes with the projections.

The upshot of this abstract nonsense (which never mentioned elements of any sets, just objects and morphisms) is that there are \emph{unique mutually inverse bijections} between $P_1$ and $P_2$ which commute with the projections to the factors $X$ and $Y$. In particular, if $P$ is \emph{a product} of $X$ and $Y$ then it is uniquely isomorphic to \emph{the product} $X\times Y$ of $X$ and $Y$ in a way compatible with the projections. In our one element example $X=\{37\}$ and $Y=\{42\}$, if $P$ is any one element set then the unique map from $P$ to $X\times Y$ of course sends the element to $(37,42)$.

\section{Products in practice.}

When a mathematician writes $X\times Y$, what do they mean? Is it \emph{a product} in the sense of the universal property, or is it the ``special'' one $X\times Y$ consisting of ordered pairs? One might imagine that, to fix our ideas, it's easiest to just choose the special one. On the other hand, a mathematician would almost certainly agree with the below claim
$$\R^2\times\R=\R\times\R^2=\R^3;$$
it is as clear as the claim that $2+1=1+2=3$. However it seems to be \emph{impossible} to set up the foundations of mathematics in such a way that all of these sets are \emph{literally equal}. Using the model of products in the previous section, a typical element of $\R^2\times\R$ looks like $((a,b),c)$ and a typical element of $\R\times\R^2$ looks like $(a,(b,c))$. These two constructions clearly carry the same data, and yet equally clearly they are not identical; they are both different models for $\R^3$, as is the model consisting of ordered triples $(a,b,c)$ defined for example as functions $\{1,2,3\}\to\R$. In particular, sets equipped with \emph{the} product do not strictly speaking form a monoid (because $(A\times B)\times C=A\times(B\times C)$ is strictly speaking false).

However all three of $\R^2\times\R$, $\R\times\R^2$ and $\R^3$ satisfy the universal property for a product of three copies of $\R$, meaning that there are \emph{unique isomorphisms} between these constructions. The category theorists would tell us that the category of sets equipped with \emph{the} product can be made into a \emph{monoidal category}, which means that we can write down the extra data of a collection of isomorphisms $i_{ABC}:(A\times B)\times C\cong A\times (B\times C)$ satisfying an equation called the \emph{pentagon axiom} \cite{wikipedia_pentagon}, which says that the two resulting natural ways of identifying $((A\times B)\times C)\times D$ with $A\times(B\times(C\times D))$ are equal. Unsurprisingly, in this example, both of the natural identifications send $(((a,b),c),d)$ to $(a,(b,(c,d)))$.

It is axioms like the pentagon axiom -- ``higher compatibitilies'' between identifications of objects which mathematicians are prone to regard as equal anyway -- which are so easy to forget. Which of $((A\times B)\times C)\times D$ and $A\times (B\times (C\times D))$ does a mathematician mean when they write $A\times B\times C\times D$? If one (strictly speaking, incorrectly) decides that the sets $((A\times B)\times C)\times D$ and $A\times(B\times(C\times D)))$ are \emph{equal} it doesn't matter! There is only one way in which two sets can be equal (in contrast to there being many ways of being isomorphic, in general), and if we think this way then we deduce the pentagon axiom no longer needs to be checked! It is phenomena like this which gives rise to arguments which are strictly speaking incomplete, throughout the literature. Note of course that in \emph{every case} known to the author, these arguments can be filled in; however the Lean community has only just started on algebraic geometry, and it will be interesting to see what happens as we progress.

I have mentioned the real numbers already. They are unique up to unique isomorphism, and mathematicians do a very good job of sticking to the universal property and developing calculus using only the completeness property of the reals rather than relying on any kind of explicit set-theoretic definition. When it comes to products however, we don't to this. Consider for example $\phi:\R^2\to\R$ defined by $\phi(x,y)=y^2+xy-x$. Mathematicians would have no objection to that definition -- however it assumes the ordered pair model for the reals: it is a function from \emph{the product} rather than from \emph{a product}. If $(P,\pi_1,\pi_2)$ is \emph{a product} then we can define $\phi_P$ on $P$ by $\phi_P(t)=\pi_2(t)^2+\pi_1(t)\pi_2(t)-\pi_1(t)$. This looks rather more ungainly than the definition of $\phi$ above so is typically avoided. However, if one wants to identify sets like $(A\times B)\times C$ and $A\times(B\times C)$ on the basis that there is a unique isomorphism between them satisfying various basic properties, then one is strictly speaking forced to develop a theory of products of sets using \emph{only} the universal property.

\section{Universal properties in algebraic geometry}\label{five}

In the previous section we saw that a simple construction in mathematics (product of two sets) could be defined up to unique isomorphism by a universal property (``a product''), and could also be defined via an explicit model (``the product''). In a perfectly functorial world, mathematicians would only ever use the universal property to do everything. However in practice this is far from the case; we are sometimes \emph{forced} to use an explicit model of our object. This surprised me. Here is an example of this phenomenon which I learnt from Patrick Massot.

Say $R$ is a commutative ring and $S$ is a multiplicatively closed subset (so $1\in S$ and $a,b\in S\implies ab\in S$). Informally, the localisation $R[1/S]$ of $R$ at $S$ is the ring generated by $R$ and the inverses of the elements of $S$. Formally there is a universal property at play here. Let us define \emph{a localisation} of $R$ at $S$ to be an $R$-algebra $R[1/S]$ with the following universal property: given any commutative ring $Z$ and any ring homomorphism $R\to Z$ sending each element of $S$ to an invertible element of $Z$, there is a unique extension of this homomorphism to an $R$-algebra map $R[1/S]\to Z$.

A localisation has a universal property, and so the uniqueness yoga for universal properties above tells us that localisations are unique up to unique isomorphism -- if they exist. To work with localisations then, it suffices to find one construction of an object with the universal property, and we could call this construction \emph{the} localisation.

Like the real numbers, there are several ways to construct an explicit model which satisfies the universal property. One is to adjoin formal inverses $y_s$ of every $s\in S$ to $R$, making a huge multi-variable polynomial ring, and then to quotient out by the ideal generated by the $sy_s-1$ for $s\in S$. Another is to consider an element of $R[1/S]$ as a fraction with numerator in $R$ and denominator in $S$, defining $R[1/S]$ to be ordered pairs $(r,s)$ modulo an appropriate equivalence relation (we say $(r_1,s_1)\sim (r_2,s_2)$ if there exists $t\in S$ such that $r_1s_2t=r_2s_1t$) and then to define addition and multiplication on the equivalence classes and to verify the commutative ring axioms. In this sense it's just like the real numbers -- we can construct them using Cauchy sequences or Dedekind cuts or Bourbaki uniform space completions, and the details don't matter. So surely here it is the same thing: now we know the \emph{existence} of a localisation, the theory should then be developed using only the universal property and we can forget about the construction. But this turns out to be essentially impossible to do in practice. Let me explain an explicit example.

It is well-known that the map $R\to R[1/S]$ making $R[1/S]$ into an $R$-algebra has kernel equal to the annihilator ideal of $S$, that is the elements $r\in R$ such that there exists $s\in S$ with $rs=0$. How does one prove this \emph{only} from the universal property of localisations? Certainly $rs=0$ implies that $r=0$ in $R[1/S]$, because $0=rss^{-1}=r1=r$ in $R[1/S]$. The harder part is the converse inclusion. Say $r\in R$ and its image in $R[1/S]$ is zero. How do we find some $s\in S$ such that $rs=0$, using just the universal property? Even the following special case seems hard: how does one prove using only the universal property that if $s\in R$ is not a zero-divisor then the map $R\to R[1/s]$ is injective?

One solution to this problem is to solve the contrapositive question. If there is no $s\in S$ such that $rs=0$, how do we find a cleverly-chosen ring $Z$ equipped with a map from $R$ to $Z$ sending $S$ to units and such that $r$ is sent to something nonzero? The way this is done in the books is as follows: let $Z$ be the model of $R[1/S]$ consisting of equivalence relations of pairs $(r,s)\in R\times S$ which we used to prove that localisations exist. Now $r$ is sent to the class of $(r,1)$, and by definition this is equivalent to the class of $(0,1)$ if and only if there exists $t$ such that $r1t=01t$, i.e. such that $rt=0$, as required.

Another way of thinking about this proof is that we have computed the kernel of \emph{the} localisation map from $R$ to \emph{the} localisation $R[1/S]$ (that is, to our explicit model of the localisation as a quotient of $R\times S$), and have then we have implicitly used an \emph{unwritten} (but easy) diagram chase to deduce that this kernel is equal the kernel of the map to \emph{a localisation}. In particular, we have violated assertion 1, because we have resorted to an explicit construction of the localisation; this is somewhat like choosing a basis for a vector space in order to reduce a question about linear maps to a question about matrices. Thanks to a referee for pointing out several other examples where it is unreasonable to expect that the universal property should do all the work. For example the claim that if $R$ is an integral domain then so is the polynomial ring $R[X]$ does not appear to immediately follow from the usual universal property describing maps $R[X]\to Z$ for other commutative rings $Z$; here it's much more natural to use the usual concrete model of $R[X]$ (functions from the naturals to $R$ with finite support, with notation $r_0+r_1X+\cdots +r_nX^n$).

Whilst picking a basis is something which sometimes cannot be avoided, my current understanding of the localisation example above is that it is evidence that we have made a \emph{poor choice} of universal property. Later on we will see another property which also uniquely characterises localisations up to unique isomorphism, and Lean's mathematics library {\tt mathlib} (which uses this choice) can be regarded as evidence that the second property is a much better choice if we want to develop the theory of localisations via universal properties and without ever choosing a model.

An unwritten but easy diagram chase might not look like a big deal, but here is an example which showed up in the Lean formalisation of algebraic geometry where this sort of thing really was an issue. In this example we only want to invert one element of a ring at a time, so let us introduce the notation $R[1/f]$ for $f\in R$ to mean $R[1/S]$ where $S=\{1,f,f^2,f^3,\ldots\}.$

In 2017 Kenny Lau (then a 1st year undergraduate at Imperial College London) constructed \emph{the} localisation of a commutative ring at a multiplicative subset in Lean (that is, he wrote down an explicit model for the localisation and proved that it satisfied the universal property). In 2018 Chris Hughes (also a 1st year undergraduate) proved \cite[\href{https://stacks.math.columbia.edu/tag/00EJ}{Tag 00EJ}]{stacks-project} in Lean. This lemma claims that if $R$ is a commutative ring and if $f_1,f_2,\ldots,f_n\in R$ generate the unit ideal then the following sequence is exact:
$$0\to R\to \prod_i R_i\to\prod_{i,j}R_{i,j}.$$
Here $i$ and $j$ (independently) run from $1$ to $n$, $R_i:=R[1/f_i]$ and $R_{i,j}:= R[1/f_if_j].$ The first nontrivial map is the structure map $R\to R_i$ on each component, and the second sends $(r_k)_k$ to $r_i-r_j$ on $R_{i,j}$.

I then tried to apply this lemma to give a formal proof in Lean that the structure presheaf on an affine scheme is a sheaf (which is the content of the assertion that an affine scheme is a scheme). However when I came to apply Hughes' lemma it \emph{would not apply} because it only applied to \emph{the localisation} of a ring at a multiplicative subset, whereas I wanted to apply it to \emph{a localisation}. More precisely, The result Hughes had formalised was about localisations such as $R[1/fg]$ and I wanted to apply it to rings such as $R[1/f][1/g]$, which are \emph{a localisation} of $R$ at the submonoid generated by $fg$, but are unfortunately not equal to \emph{the localisation}. For example if we use the $R\times S$ model for \emph{the localisation}, an element of $R[1/fg]$ is an equivalence class in $R\times\{1,fg,(fg)^2,\ldots\}$, and an element of $R[1/f][1/g]$ is an equivalence class in the completely different set $R[1/f]\times\{1,\overline{g},\overline{g}^2,\ldots\}$, where $\overline{g}$ is the image of $g$ in $R[1/f]$ (so itself is an equivalence class). 

The bottom line is that the lemma did not apply. And yet mathematically there should be no issue at all. How do we fix this problem?

Let me go over the issue again. The statement of the ring lemma which Hughes had formalised was a theorem about \emph{the localisation} of a ring at many different multiplicative sets, but in our application we needed the stronger statement that it was true for \emph{a localisation}. Mathematicians barely notice the difference, but a computer theorem prover was quick to point it out.

The cheap solution would be to \emph{deduce} the stronger statement from the weaker one, and this is what I did at first. This sounds like it should be a relatively straightforward process, and it is arguably what is happening inside the head of a mathematician who is showing this argument to their algebraic geometry class; it involves introducing no new mathematical ideas, it is just a big diagram chase. However, to my dismay, this was in practice \emph{extremely} tedious. I had to prove lemma after tedious lemma saying that various diagrams commuted, using the universal property; these lemmas were specifically about morphisms which appeared in the statement of the ring lemma, and in particular were not reusable in other situations; they were just annoying boilerplate which needed to be in the code base and furthermore they felt mathematically completely vacuous. However I could see that they were necessary to make this proof strategy watertight. Furthermore we ran into an unexpected twist: the map from $\prod_k R_k$ to $\prod_{i,j}R_{i,j}$ in the lemma sending $(r_k)$ to $r_i-r_j$ is not an $R$-algebra homomorphism but only an $R$-module homomorphism, and there can be many (non-canonical) $R$-module maps between two localisations of $R$. It is however the difference of two $R$-algebra homomorphisms, and applying the universal property to these ring homomorphisms got us through in the end.

The major disadvantage of this approach is that the human formaliser has to generate a bunch of one-off lemmas with no applications elsewhere. It took a while for me to understand that in fact formalising the lemma as stated in the literature, for \emph{the localisations}, is a mistake. One should in fact never prove this lemma at all; one should instead state and prove a version of the lemma which assumes only that the $R_i$ and $R_{i,j}$ are each \emph{a} localisation of $R$. If one likes, one can then deduce the weaker version (which is of no use to us anyway) from the stronger one. In 2018 Amelia Livingston, also then an undergraduate at Imperial College, developed theories of both ``the'' and ``a'' localisation of commutative monoids at submonoids, and then a theory of ``the'' and ``a'' localisation of commutative rings at submonoids (thus refactoring Lau's work). Finally in 2019 Ramon Fern\'andez Mir, as part of his MSc thesis, came up with a proof of the version of the ring lemma which applied to rings which were merely ``a'' localisation, and we could use these results to give a new and much cleaner Lean proof that the structure presheaf on an affine scheme was a sheaf. Mir followed a suggestion of Neil Strickland to instead use a \emph{different} property which also uniquely characterises localisations, and then to find a proof of the ring theory lemma which used only this new property. The property Mir used was the following:

\begin{theorem*} The $R$-algebra $A$ with structure homomorphism $f:R\to A$ is a localisation of $R$ at a multiplicative subset $S$ if and only if it satisfies the following three properties:
  \begin{itemize}
  \item $f(s)$ is invertible for all $s\in S$;
  \item Every element of $A$ can be written as $f(r)/f(s)$ for some $r\in R$ and $s\in S$;
  \item The kernel of $f$ is the annihilator of $S$.
  \end{itemize}
\end{theorem*}

A referee notes that there is an analogous classification of localisations of $R$-modules, which we leave as an exercise for the reader. This property is of a very different nature to the previously-mentioned universal property; for example it does not quantify over all rings, meaning that when applying this theorem to prove results about localisations, one never has to pull a magic ring out of a hat (such as an explicit construction of \emph{the} localisation). It does however uniquely characterise the localisation up to unique isomorphism, so it is doing the job required of it. It is also a mathematical idea which is \emph{not present} in the original paper presentation of the proof; although it is not difficult to prove, it is in some sense \emph{new mathematics} which had to be generated to make the formalisation viable and tidy. It also explicitly states the characterisation of the kernel of $f$ which Massot had observed seemed to be difficult to prove from the traditional universal property without resorting to an explicit model. Mir was fortunate in that the proof formalised by Hughes for \emph{the} localisation was easily adaptable so that it applied to all localisations satisfying Strickland's property; it was not at all a priori clear (at least to me) that this would be the case.

The very fact that we needed to do some mathematical thinking (the invention of a practical property characterising localisations up to unique isomorphism) in order to come up with an acceptable formalisation of this basic result in algebraic geometry made me wonder just how much more mathematical thinking we will need to do in order to formally prove harder results in algebraic geometry.

\section{The problem with Grothendieck's use of equality.}

The above story is evidence that there is a missing argument in the literature, and that the statement of \cite[\href{https://stacks.math.columbia.edu/tag/00EJ}{Tag 00EJ}]{stacks-project} in the stacks project is, strictly speaking, not strong enough deduce the claim that the structure presheaf is a sheaf just before \cite[\href{https://stacks.math.columbia.edu/tag/01HU}{Tag 01HU}]{stacks-project}. However the algebraic geometry community does not regard this issue as problematic, and indeed the way the theory is presented it is extremely difficult to even notice that this is an issue. 
I believe that one major reason for this can be traced back to Grothendieck's seminal work EGA (\cite{EGA1}), where he and Dieudonne develop the foundations of modern algebraic geometry. The word ``canonique'' appears \emph{hundreds} of times in EGA1, with no definition ever supplied. The argument which caused all our trouble is in section 1.3. Grothendieck claims that if $R$ is a commutative ring and $f,g$ are two elements contained in the same prime ideals (for example we could have $R=\C[X]$, $f=X^2$ and $g=X^3$) and if $S$ is the multiplicative subset of $R$ containing the $s$ which divide some power $f^n$ of $f$ (or equivalently divide some power of $g$) then $R[1/f]$ and $R[1/g]$ both ``s'identifient canoniquement'' with ring $R[1/S]$ and hence $R[1/f]=R[1/g]$ (see for example section 1.3.3 of~\cite{EGA1} where the stronger statement $M[1/f]=M[1/g]$ is claimed for any $R$-module $M$). Of course we certainly know what Grothendieck \emph{means} -- $R[1/f]$ and $R[1/g]$ are uniquely isomorphic as $R$-algebras, so we will identify them via this isomorphism and then call it an equality. Lean would tell Grothendieck that this equality \emph{simply isn't true} and would stubbornly point out any place where it was used. Let me emphasize once more: Grothendieck was well aware of what he was saying, but Lean would argue that he was confusing $=$ and $\cong$.

The idea that objects could be ``canonically'' isomorphic seems to have been taken on with some enthusiasm by many in the mathematical community\footnote{Perhaps one exception was Andr\'e Weil: in a letter to the editor of the Annals\cite{hades} (purporting to be from Lipschitz writing from Hades), he says, with a clear degree of sarcasm, ``I can assure you, at any rate, that my intentions are honourable and my results invariant, probably canonical, perhaps even functorial.''}. By the 1970s it was clear that Grothendieck's ideas were here to stay: his discovery of \'etale cohomology had led to a proof of the Weil conjectures, fundamental statements about the number of solutions to polynomial equations over finite fields which could be made without any reference to the theory of schemes, but which apparently could only be proved using them. The book \cite{Milne} from 1980, one of the first textbook treatments of etale cohomology, contains in its ``Terminology and conventions'' section, the convention that ``a canonical isomorphism [is denoted by] $=$''. Nowhere in any of these texts is any definition of the word ``canonical''. Gordon James told me that he once asked John Conway what the word meant, and Conway's reply was that if you and the person in the office next to yours both write down a map from $A$ to $B$, and it's the same map, then this map is canonical. This might be a good joke, but it is not a good definition.

In Milne's book we do not just see localisations -- we see pullbacks and pushforwards of sheaves, maps defined as adjoint functors, we see limits, colimits, quotients by equivalence relations, tensor products of modules, and constructions coming from Grothendieck's six functor formalism. All of these constructions are universal and no doubt any maps produced by these universal properties are ``canonical''. But whenever Milne is (ab)using the equality symbol, there should \emph{in theory} be a check that whatever theory is being developed is valid for any object satisfying the universal property in question. The reason that this is not happening is the devious technique of arguing that two objects which satisfy the same universal property are \emph{``canonically'' isomorphic} and hence ``can be identified'' and hence ``are equal''. To give a random example from~\cite{Milne}: in Section II.3, Remark 3.1(f) Milne talks about the direct and inverse image of a sheaf under a morphism of sites, and claims that $(\pi'\pi)_*=\pi'_*\pi_*$ and $\pi^*\pi'^*=(\pi'\pi)^*.$ Equality of functors is in some sense not a sensible mathematical notion, as it boils down to many statements about equality of objects in a category, and equality of objects is not invariant under equivalence of categories -- it is hence sometimes referred to as an ``evil'' concept for this reason. Moreover, as well as being evil, the claim is \emph{not actually correct} for pullbacks, because ``the'' pullback of a sheaf involves making a choice of an explicit construction of sheafification of a presheaf, and the set-theoretic equality $\pi^*\pi'^*F=(\pi'\pi)^*F$ fails for essentially the same reason as the equality $R[1/f][1/g]=R[1/fg]$ fails. What is actually going on is a functorial identification which satisfies some unwritten compatibilities -- the details are ``left to the reader''. I thank the referee for pointing out that the implicit use of natural isomorphisms in sheaf cohomology when dealing with pullback functors goes back to Godement's book on sheaf theory. Before one formalises this kind of mathematics, one will have to think carefully about precisely which of the properties which uniquely define pullbacks of sheaves up to unique isomorphism should be used as the definition of \emph{a} pullback, and, just as in the case of localisation of rings, it might not be the one which first springs to mind: the definition of pullback as being adjoint to pushforward involves quantification over all sheaves on a site and hence might not be the most ergonomic characterisation.

\section{More on ``canonical'' maps}

The previous remarks have been mostly the flagging of a technical point involving mathematicians ``cheating'' by considering that various nonequal but uniquely isomorphic things are equal, and a theorem prover pointing out the gap. Whilst I find this subtlety interesting, I do not believe that this slightly dangerous convention is actually hiding any errors in algebraic geometry; all it means is that in practice people wishing to formalise algebraic geometry in theorem provers are going to have to do some work thinking hard about universal properties, and possibly generate some new mathematics in order to make the formalisation of modern algebraic geometry a manageable task. Section 1.2 of Conrad's book~\cite{conrad_duality} gives me hope; his variant of the convention is summarised there by the following remark: ``We sometimes write $A=B$ to denote the fact that $A$ is canonically isomorphic to $B$ (via an isomorphism which is always clear from the context).'' Even though we still do not have a definition of ``canonical'', we are assured that, throughout Conrad's work at least, it will be clear which identification is being talked about.

In the work of Grothendieck we highlighted, the rings he calls ``canonically isomorphic'' are in fact \emph{uniquely} isomorphic as $R$-algebras. However when it comes to the Langlands Program, ``mission creep'' for the word ``canonical'' is beginning to take over. Before I discuss an example from the literature let me talk about a far more innocuous use of the word. Consider the following claim:

\begin{theorem*}[The first isomorphism theorem] If $\phi:G\to H$ is a group homomorphism, then $G/\ker(\phi)$ and $\im(\phi)$ are canonically isomorphic.
\end{theorem*}

I think that we would all agree that the first isomorphism theorem does say strictly more than the claim that $G/\ker(\phi)$ and $\im(\phi)$ are \emph{isomorphic} -- the theorem is attempting to make the stronger claim that there is a ``special'' map from one group to the other (namely the one sending $g\ker(\phi)$ to $\phi(g)$) and that it is \emph{this} map which is an isomorphism. In fact this is the claim which is used in practice when applying the first isomorphism theorem -- the mere existence of an isomorphism is often not enough; we need the formula for it. We conclude

\begin{theorem*} The first isomorphism ``theorem'' as stated above is not a theorem.\end{theorem*}

Indeed, the first isomorphism ``theorem'' is a \emph{pair} consisting of the \emph{definition} of a group homomorphism $c:G/\ker(\phi)\to\im(\phi)$, and a \emph{proof} that $c$ is an isomorphism of groups. In contrast to earlier sections, uniqueness of the isomorphism is now \emph{not true} in general. For example, if $H$ is abelian, then the map $c^*$ sending $\overline{g}\in G/\ker(\phi)$ to $c(\overline{g})^{-1}$ is also an isomorphism of groups, however this isomorphism is not ``canonical'': an informal reason for this might be ``because it contains a spurious ${}^{-1}$'', but here a better reason would be because it does not commute with the canonical maps from $G$ to $G/\ker(\phi)$ and $H$. What is actually going on here is an implicit \emph{construction}, as well as a theorem. The claim implicit in the ``theorem'' is that we can \emph{write down a formula} for the isomorphism -- we have \emph{made} it, rather than just deduced its existence from a nonconstructive mathematical fact such as the axiom of choice or the law of the excluded middle. My belief is that some mathematicians have lost sight of this point, and hence are confusing constructions (definitions) with claims of ``canonical''ness (attempts to state theorems). The currency of the mathematician is the theorem, so theorems we will state.

\section{Canonical isomorphisms in more advanced mathematics}

We finish this paper by taking a look at some more advanced topics. We start with the Langlands program. In Langlands' 1968 paper~\cite{langlands_abelian} he proves a theorem which is now referred to as ``the proof of the local and global Langlands conjectures for abelian algebraic groups'', although as Langlands himself explains, the result is an exercise using standard results from local and global class field theory. The set-up is that $K/F$ is a finite Galois extension of fields, $T$ is a torus over $F$ which splits over $K$, $L$ is the character group of~$T$ with its action of $\Gal(K/F)$, and $\widehat{T}$ is the dual complex torus of~$T$ equipped with the induced action of $\Gal(K/F)$. We directly quote part of the first theorem in the paper.

\medskip

``{\bf Theorem 1.} If $K$ is a global or local field there is a canonical isomorphism of $H^1_c(W_{K/F},\widehat{T})$ with the group of generalized characters of $\Hom_{\Gal(K/F)}(L,C_K)$.''

\medskip

Here, if $K$ is local or global, then the Weil group $W_{K/F}$ is a certain topological group equipped with a map to $\Gal(K/F)$, and $H^1_c$ denotes continuous cocycles over coboundaries.

The issue with this ``theorem'' is that, just like the first isomorphism theorem, it is \emph{not a theorem} -- it is more than this. This theorem cannot be applied without reading the proof: this is in fact a sure fire sign that it is not actually a theorem. From a formalist point of view, theorems are mathematical objects which are completely encapsulated by their \emph{statements}. Here is an argument to show that there must be more going on in this ``theorem''. The claim is that one group (a cohomology group) is canonically isomorphic to another group (a character group). Both groups in question are abelian, and are coming from completely different worlds; the cohomology group $H^1_c(W_{K/F},\widehat{T})$ is a Galois-theoretic object and hence ``algebraic''; the generalized characters are an automorphic object and hence ``analytic''. What Langlands gives in the proof, of course, is the explicit \emph{construction} of a map from the cohomology group to the character group; part of his contribution is a \emph{definition} of this map via an ``explicit formula''; it is not completely constructive (in the sense of constructive mathematics) because it involves the inverse of a map which is only shown to be a bijection via a nonconstructive cohomological argument, and the inverse of a constructive bijection might not be constructive (see~\cite{xena-bijection}). However, it is an explicit \emph{definition}. Let's call it $d$. Now consider the function $d^*$ obtained by composing $d$ with the inversion map on the target abelian group. This is still a bijection. Is it also ``canonical''? This question cannot be answered, because the word ``canonical'' has no formal mathematical definition. However, unlike the example of the first isomorphism theorem, here the problem is genuinely worse.

The issue in fact goes back to class field theory, where one of the main theorems is a ``canonical'' isomorphism between an abelianised global Galois group and an adelic group. Again both of these topological groups are abelian, so given one ``canonical'' isomorphism we can compose with inversion and get another one. Both of these ``canonical'' isomorphisms are known to mathematicians -- indeed they have different names! One is them is ``the isomorphism of class field theory sending an arithmetic Frobenius element to a uniformiser'' and the other is ``the isomorphism of class field theory sending a geometric Frobenius element to a uniformiser''; a ``geometric Frobenius element'' is just defined to be the inverse of an arithmetic Frobenius element. Which of these maps is the ``canonical'' one? In my experience, both of these normalisations are used. For example people working in the area of the cohomology of Shimura varieties use the geometric normalisation, but people working in the theory of Heegner points (and hence with Tate modules, which are homology groups) use the arithmetic one, in both cases because it locally minimises the number of minus signs in the results. There seems to be no good answer to the question of which isomorphism is ``preferred'' in mathematics -- different mathematicians working in different areas prefer different choices, and have coherent reasons to prefer their choice over the other.

Another place where sign issues pervade mathematics is in the theory of homological algebra, also now a key tool being used in the Langlands program. Let us consider a very simple case, namely group cohomology. If $G$ is a group and we have a short exact sequence $0\to A\to B\to C\to 0$ of $G$-modules, what is ``the'' definition of ``the'' boundary map $H^0(G,C)\to H^1(G,A)$? One starts with a $G$-invariant element $c\in C$, lifts it to $b\in B$, and now given $g\in G$ we observe that $b$ and $gb$ both map to $c$, so their difference maps to zero and is hence in $A$. The 1-cocycle representing the cohomology class sends $g$ to this difference. But is this difference $b-gb$ or $gb-b$? Which is the ``canonical'' choice? It has taken me a long time to realise that this question \emph{does not have a preferred answer}. Grothendieck was well aware of this; his concept of a universal $\delta$-functor solves this problem. As \emph{part of the data} of a universal $\delta$-functor one has to \emph{supply} the boundary homomorphisms (in all degrees); they do not come to us by magic. As Grothendieck knew, it is not enough to define a cohomology theory by defining the cohomology groups; the extra data of the boundary homomorphisms must also be supplied: these are only ``canonical up to a unit'', and $-1$ is a unit in the integers.

Signs are also a nightmare whenever one is collapsing a double complex to a single complex; there is no ``canonical'' way to do this, and a choice must be made. Is the same choice made throughout the literature? The answer, unfortunately, is ``no''. This means that extreme care must be taken when quoting results from more than one reference in this area. Again I refer the reader to Conrad's book~\cite{conrad_duality}, where in the introduction great lengths are taken to explain incompatibilities in the references he wants to cite. Strictly speaking, this is how all mathematics should look, unpleasant as it seems. And when it comes to formalisation of this mathematics, as it one day will, these things really will matter.

\section{Summary}

Whilst I am not making any claims about errors in the literature or holes in arguments which cannot be filled in after some work, I am arguing that these holes do exist, and that new mathematics might need to be done by formalisers in order to fill in these holes in an efficient way.

The holes are of two kinds. Firstly there is the issue of people making constructions or proving theorems which make essential use of a model of a mathematical object which is defined up to unique isomorphism; the hole here is that it needs to be checked that the argument does not depend on the explicit details of the model. Mathematicians are well aware of this when it comes to, say, picking a basis for a vector space and then checking that nothing important depended on the choice, or picking a representative for an equivalence class and then checking that nothing important depends on the representative. However they seem to be less careful when doing more advanced mathematics, confusing ``a'' localisation with ``the'' localisation or ``a'' pullback with ``the'' pullback, and leaving to the reader the details of checking that many diagrams commute. One useful trick is to abuse the equality symbol, making it mean something which it does not mean; this can trick the reader into thinking that nothing needs to be checked. Sometimes such checks can be surprisingly painful, and it may be easier to restructure a mathematical argument than to actually make these checks.

The second kind of hole is the issue of various maps (like boundary maps in exact sequences) being regarded as ``canonical'' where now they are in fact \emph{not} unique, and there are implicit choices of sign being made. Unfortunately it is not at all difficult to point to explicit examples in the literature where an author does not state precisely which convention they are using when it comes to things like the theory of Shimura varieties, or homological algebra. This puts an unnecessary burden on the careful mathematician (for example Conrad, or a computer theorem prover) who is attempting to use or verify the work.

Both of these issues have shown up in my formalisation work, and I expect them to show up more often as we go deeper into the formalisation of modern mathematics.

\bibliographystyle{amsalpha}

\bibliography{equal.bib}

\providecommand{\bysame}{\leavevmode\hbox to3em{\hrulefill}\thinspace}
\providecommand{\MR}{\relax\ifhmode\unskip\space\fi MR }
\providecommand{\MRhref}[2]{%
  \href{http://www.ams.org/mathscinet-getitem?mr=#1}{#2}
}
\providecommand{\href}[2]{#2}
\begin{thebibliography}{{Wik}04b}

\bibitem[AX23]{angdinata_yu}
David~Kurniadi Angdinata and Junyan Xu, \emph{{An Elementary Formal Proof of
  the Group Law on Weierstrass Elliptic Curves in Any Characteristic}}, 14th
  International Conference on Interactive Theorem Proving (ITP 2023) (Dagstuhl,
  Germany) (Adam Naumowicz and Ren\'{e} Thiemann, eds.), Leibniz International
  Proceedings in Informatics (LIPIcs), vol. 268, Schloss Dagstuhl --
  Leibniz-Zentrum f{\"u}r Informatik, 2023, pp.~6:1--6:19.

\bibitem[BCM20]{bcmperfectoid}
Kevin Buzzard, Johan Commelin, and Patrick Massot, \emph{Formalising perfectoid
  spaces}, Proceedings of the 9th {ACM} {SIGPLAN} International Conference on
  Certified Programs and Proofs, {CPP} 2020, New Orleans, LA, USA, January
  20-21, 2020 (Jasmin Blanchette and Catalin Hritcu, eds.), {ACM}, 2020,
  pp.~299--312.

\bibitem[BPL21]{isabelle-schemes}
Anthony Bordg, Lawrence Paulson, and Wenda Li, \emph{Grothendieck's schemes in
  algebraic geometry}, March 2021,
  \url{https://isa-afp.org/entries/Grothendieck_Schemes.html}, Formal proof
  development.

\bibitem[Buz]{chapmanvid}
Kevin~M. Buzzard, \emph{Grothendieck's approach to equality},
  \url{https://www.youtube.com/watch?v=-OjCMsqZ9ww}, Accessed: 12-08-2023.

\bibitem[{Buz}19]{xena-bijection}
{Buzzard, Kevin}, \emph{The inverse of a bijection}, 2019, [Online; accessed
  12-Aug-2023].

\bibitem[Con00]{conrad_duality}
Brian Conrad, \emph{Grothendieck duality and base change}, Lecture Notes in
  Mathematics, vol. 1750, Springer-Verlag, Berlin, 2000. \MR{1804902}

\bibitem[dFF23]{defrutosfernandeznorms}
Mar{\'\i}a~In\'{e}s de~Frutos-Fern\'{a}ndez, \emph{{Formalizing Norm Extensions
  and Applications to Number Theory}}, 14th International Conference on
  Interactive Theorem Proving (ITP 2023) (Dagstuhl, Germany) (Adam Naumowicz
  and Ren\'{e} Thiemann, eds.), Leibniz International Proceedings in
  Informatics (LIPIcs), vol. 268, Schloss Dagstuhl -- Leibniz-Zentrum f{\"u}r
  Informatik, 2023, pp.~13:1--13:18.

\bibitem[Gro60]{EGA1}
A.~Grothendieck, \emph{\'{E}l\'{e}ments de g\'{e}om\'{e}trie alg\'{e}brique.
  {I}. {L}e langage des sch\'{e}mas}, Inst. Hautes \'{E}tudes Sci. Publ. Math.
  (1960), no.~4, 228. \MR{217083}

\bibitem[Lan97]{langlands_abelian}
R.~P. Langlands, \emph{Representations of abelian algebraic groups}, Pacific J.
  Math. (1997), 231--250, Olga Taussky-Todd: in memoriam. \MR{1610871}

\bibitem[Liv23]{livingstongpcoh}
Amelia Livingston, \emph{{Group Cohomology in the Lean Community Library}},
  14th International Conference on Interactive Theorem Proving (ITP 2023)
  (Dagstuhl, Germany) (Adam Naumowicz and Ren\'{e} Thiemann, eds.), Leibniz
  International Proceedings in Informatics (LIPIcs), vol. 268, Schloss Dagstuhl
  -- Leibniz-Zentrum f{\"u}r Informatik, 2023, pp.~22:1--22:17.

\bibitem[mC20]{mathlib}
The mathlib Community, \emph{The lean mathematical library}, Proceedings of the
  9th {ACM} {SIGPLAN} International Conference on Certified Programs and
  Proofs, {ACM}, jan 2020.

\bibitem[Mil80]{Milne}
James~S. Milne, \emph{\'{E}tale cohomology}, Princeton Mathematical Series, No.
  33, Princeton University Press, Princeton, N.J., 1980. \MR{559531}

\bibitem[{Sta}18]{stacks-project}
The {Stacks Project Authors}, \emph{\textit{Stacks Project}},
  \url{https://stacks.math.columbia.edu}, 2018.

\bibitem[Wei59]{hades}
André Weil, \emph{Correspondence [signed ``{R}. {L}ipschitz'']}, Ann. of Math.
  (2) \textbf{69} (1959), 247--251, Attributed to A. Weil. \MR{100637}

\bibitem[{Wik}04a]{wikipedia_pentagon}
{Wikipedia contributors}, \emph{Monoidal category --- {W}ikipedia{,} the free
  encyclopedia}, 2004, [Online; accessed 12-Aug-2023].

\bibitem[{Wik}04b]{wikipedia_ordered_pair}
\bysame, \emph{Ordered pair --- {W}ikipedia{,} the free encyclopedia}, 2004,
  [Online; accessed 20-May-2023].

\bibitem[Zha23]{zhang_proj}
Jujian Zhang, \emph{{Formalising the Proj Construction in Lean}}, 14th
  International Conference on Interactive Theorem Proving (ITP 2023) (Dagstuhl,
  Germany) (Adam Naumowicz and Ren\'{e} Thiemann, eds.), Leibniz International
  Proceedings in Informatics (LIPIcs), vol. 268, Schloss Dagstuhl --
  Leibniz-Zentrum f{\"u}r Informatik, 2023, pp.~35:1--35:17.

\bibitem[ZM23]{zeuner2023univalent}
Max Zeuner and Anders Mörtberg, \emph{A univalent formalization of
  constructive affine schemes}, 2023.

\end{thebibliography}

\end{document}